\documentclass[12pt]{amsart}

\usepackage{amsmath,amssymb}
\usepackage{accents}
\usepackage{graphicx}
\usepackage{verbatim}
\usepackage{color}
\usepackage{float}
\usepackage{enumerate}

\newlength{\dhatheight}

\newtheorem{theorem}{Theorem}[section]
\newtheorem{conjecture}[theorem]{Conjecture}

\newtheorem{corollary}[theorem]{Corollary}
\newtheorem{question}[theorem]{Question}

\begin{document}
\title{Unknotting number is not additive under connected sum}

\author{Mark Brittenham and Susan Hermiller}
\address{Department of Mathematics\\
        University of Nebraska\\
         Lincoln NE 68588-0130, USA\\
        mbrittenham2@unl.edu}
\address{Department of Mathematics\\
        University of Nebraska\\
         Lincoln NE 68588-0130, USA\\
        hermiller@unl.edu}

\date{September 13, 2025}

\maketitle

\begin{abstract}
We give the first examples of a pair of knots $K_1$,$K_2$ in the 3-sphere 
for which their unknotting numbers satisfy
$u(K_1\#K_2)<u(K_1)+u(K_2)$ . This answers question 1.69(B) from Kirby's 
problem list \cite{kir93}, in the negative.
\end{abstract}

%%\keywords{unknotting number, connected sum}


%%%%%%%%%%%%%%%%%%%%%%%%%%%%%%%%%%%%%%%%%%%%%%%%%%%%%%%%%%%%%%%%%%%%%%%%%%%%
%%%%%%%%%%%%%%%%%%%%%%%%%%%%%%%%%%%%%%%%%%%%%%%%%%%%%%%%%%%%%%%%%%%%%%%%%%%%

\section{Introduction}\label{sec:intro}

%%%%%%%%%%%%%%%%%%%%%%%%%%%%%%%%%%%%%%%%%%%%%%%%%%%%%%%%%%%%%%%%%%%%%%%%%%%%
%%%%%%%%%%%%%%%%%%%%%%%%%%%%%%%%%%%%%%%%%%%%%%%%%%%%%%%%%%%%%%%%%%%%%%%%%%%%

For knots $K$ in the 3-sphere $S^3$, the \emph{unknotting number} (or \emph{Gordian number}) 
$u(K)$ of $K$ is one of the most fundamental
measures of the complexity of the knot. It is defined as the 
minimum number of crossing changes, interspersed with isotopy,
required to transform a diagram of $K$ to a diagram of the unknot. A
wide array of techniques have been employed over the years
to compute unknotting numbers, using nearly every new technique
that has been introduced into knot theory (see 
\cite{coli}, \cite{kawa}, \cite{krmr}, \cite{li02}, \cite{mc17}, \cite{mu65}, \cite{ow08},
\cite{ow10}, \cite{ozsz}, \cite{ru93}, \cite{sto02}
for a selection). Yet to this day there
are still nine prime 10-crossing knots whose unknotting numbers remain unknown,
a testament to the difficulty of determining this basic invariant.
See the Knotinfo site \cite{knotinfo} for the most up-to-date list of unresolved unknotting numbers.

Connected sum is a basic operation for creating or decomposing knots and links.
A diagram for the connected sum $K_1\#K_2$ of two knots $K_1$ and $K_2$ can be obtained from 
disjoint diagrams of the two knots by deleting a short arc in each and 
introducing two new arcs to connect the endpoints of the deleted arcs to
create a knot, without introducing additional crossings in the diagram. 
A knot is \emph{prime} if it cannot be expressed as the connected sum of 
two non-trivial knots. Every knot admits a unique decomposition (up to ordering) 
as a connected sum of prime knots \cite{schu61}. This has enabled many problems
in knot theory to focus on prime knots as their principle object of study.

Many knot-theoretic constructions and invariants are known to behave well under 
the connected sum operation. Polynomial invariants, such as the Alexander \cite{alex28}, 
Jones \cite{jones85}, and HOMFLY 
polynomials \cite{homfly} are multiplicative under connected sum, and many 
basic numerical invariants are very well-behaved under connected sum, 
such as the 3-genus \cite{schu53}, signature \cite{mu65}, braid index \cite{bm90},
and bridge number \cite{schu54}. Unknotting number has long been 
conjectured to be additive under connected sum; this conjecture is implicit 
in the work of Wendt,
in one of the first systematic studies of unknotting number \cite{wendt37}.
It is unclear when and where this was first explicitly stated; most references
to it call it an `old conjecture'. It can be found in the problem list
of Gordon \cite{plans77} from 1977 and in Kirby's list \cite{kir93}.

\smallskip

\begin{conjecture}\label{conj:additivity}
For knots $K_1$ and $K_2$ in the 3-sphere, $u(K_1\#K_2)=u(K_1)+u(K_2)$.
\end{conjecture}

The inequality $u(K_1\#K_2)\leq u(K_1)+u(K_2)$ is 
always true, and is very straightforward; choose 
diagrams for $K_1$ and $K_2$ that can be unknotted with $u(K_1)$ and
$u(K_2)$ crossing changes, respectively, and create the connected sum 
$K_1\#K_2$ using these
diagrams. Then the same set of $u(K_1)+u(K_2)$ crossing changes
will transform $K_1\#K_2$ into (unknot)$\#$(unknot) = the unknot.
The challenge has always been, essentially, to determine that one 
cannot do better
than this approach, and so $u(K_1\#K_2)\geq u(K_1)+u(K_2)$, as well.

Progress on this conjecture has been slow, however. 
Several authors have found families
${\mathcal K}$ of knots for which additivity holds for pairs of knots
chosen within the family \cite{krmr},\cite{ba05}. Perhaps the most
significant work in support of the conjecture has been Scharlemann's result 
\cite{scha85} that a knot $K$ with $u(K)=1$ is prime; consequently, 
the connected sum of two non-trivial knots must have unknotting number 
at least 2. So unknotting number is additive, when both summands have 
unknotting number 1. 
It is also well-known that if two knots $K_1$ and $K_2$ satisfy
$|\sigma(K_i)|=2u(K_i)$ for each $i$, where $\sigma(K)$ is the signature of $K$, 
and $\sigma(K_1)$ and $\sigma(K_2)$ have the same sign, then 
$u(K_1\#K_2)=u(K_1)+u(K_2)$. This is a consequence of the inequality
$|\sigma(K)|\leq 2u(K)$ \cite{mu65} and the additivity of signature under
connected sum. Since $|\sigma(K_i)|=2u(K_i)$ is actually quite common
among low crossing number knots (see, e.g., \cite{knotinfo}), this
offers strong support to the conjecture.

More recently, Yang \cite{ya08} has shown that the
connected sum of $n$ identical knots $K$ whose 
Alexander polynomial satisfies $\Delta_K(t)\neq 1$
has unknotting number at least $n$. Very recently, Applebaum, Blackwell, 
Davies, Edlich, Juh\'{a}sz, Lackenby, Toma\v{s}ev, and Zheng
\cite{abdejltz24} employed reinforcment learning to search for 
counterexamples to Conjecture \ref{conj:additivity}, and established some interesting
consequences for unknotting sequences of connected sums. They also
showed how, if Conjecture \ref{conj:additivity} were true,
one could use this to compute the unknotting numbers of many prime knots.

It was this last perspective - what can we learn if we assume that 
unknotting number is additive - which inspired the authors of the present work to carry
out the computations that led to their discovery. In the following we use the notation $\overline{K}$
to represent the mirror image of the knot $K$. It is an elementary
fact that $u(K)=u(\overline{K})$ for every knot $K$.

\begin{theorem}\label{thm:counter}
Unknotting number is not additive under connected sum.
In particular, for the knot $K=7_1$ (the (2,7)-torus knot), with
$u(K)=3$, the connected sum $L=K\#\overline{K}$ satisfies 
$u(L)\leq 5<6=u(K)+u(\overline{K})$.
\end{theorem}

The additivity of unknotting number appears as Problem 1.69(B) in 
Kirby's problem list \cite{kir93}; Theorem \ref{thm:counter} therefore settles
this question (in the negative). 

\smallskip

Once we have this one example, it takes little work to find infinitely many more.
The torus knot $7_1=T(2,7)$ appears in a minimal unknotting sequence for
each of the torus knots $T(2,2k+1)$, $k\geq 3$; see Section \ref{sec:counter}.
From this we can show:

\begin{corollary}\label{cor:torus}
For every natural number $k,\ell\geq 3$, 
$$u(T(2,2k+1)\#\overline{T(2,2\ell+1)})<u(T(2,2k+1))+u(T(2,2\ell+1))\ .$$
\end{corollary}

\noindent Section \ref{sec:counter} describes further examples.

The knot group $\pi(K)=\pi_1(S^2\setminus K)$, the fundamental group of the
complement of the knot in the $3$-sphere, is a powerful tool for the 
study of a knot. For many knots, 
it determines the knot up to isotopy and mirroring \cite{wh87},\cite{gl89}. However,
the knot group does not detect mirroring, and, in particular, 
$\pi(K\# L)\cong\pi(K\#\overline{L})$ for any knots $K$ and $L$
(see, for example, \cite[Corollary~6.1.10]{ka96}). Since
$\sigma(7_1\# 7_1)=12$, and so $u(7_1\# 7_1)=6$, while $u(7_1\#\overline{7_1})\leq 5$, 
we have

\begin{corollary} \label{cor:knotgroup}
The knot group $\pi(K)$ of a knot $K$ does not determine the unknotting number of $K$.
\end{corollary}

Using work of Baader \cite{ba06} we can show that, using the examples of Corollary 
\ref{cor:torus}, there are infinitely many non-torus knots for which the theorem holds, as well.

\begin{corollary} \label{cor:all_unk_nums}
For every $n\geq 4$ there are infinitely many knots $K_i$
with $u(K_i)=n$ and $u(K_i\#\overline{K_j})<u(K_i)+u(\overline{K_j})$ for every $i$ and $j$.
\end{corollary}

By taking connected sums of knots idenified in the above corollaries, we can make the additivity
gap arbitrarily large.

\begin{corollary} \label{cor:large_gaps}
For every $n\geq 1$ there are infinitely many knots $K_i$
with $u(K_i\#\overline{K_i})\leq u(K_i)+u(\overline{K_i})-n$.
\end{corollary}

Further methods for building knots with large gap are discussed in Section \ref{sec:road}.

Our approach to discovering our main example, described in section \ref{sec:howfound} below, is modeled 
on the authors' previous work settling the Bernhard-Jablan Unknotting Conjecture \cite{bh21}, 
and at least partly relies on data generated from that project. We 
also must credit, and express our immense gratitude to, the authors 
and developers of SnapPy/SnapPea \cite{snappy}; no
part of this project would have been possible without the extensive
computations made possible by this software. 

\medskip

The structure of this paper is as follows. In Section \ref{sec:counter} we discuss the 
knots $K=7_1$ and $L=K\#\overline{K}$, and give the knot diagrams
and crossing changes which establish the failure of additivity for unknotting number. 
We also show how to use this result to find many more examples where additivity fails.
In Section \ref{sec:howfound} we describe the motivation for our searches, 
including work on a related conjecture, on the invariance of unknotting number 
under knot mutation 
(which is Problem 1.69(C) of \cite{kir93}), and outline the process
by which we carried the searches out. In Section \ref{sec:road} we discuss 
the consequences
of this example for the computation of unknotting number, and possible 
future directions. The final section provides python code which will allow 
the reader to carry out their
own verification of the computations used in our main example.

%%%%%%%%%%%%%%%%%%%%%%%%%%%%%%%%%%%%%%%%%%%%%%%%%%%%%%%%%%%%%%%%%%%%%%%%%%%%
%%%%%%%%%%%%%%%%%%%%%%%%%%%%%%%%%%%%%%%%%%%%%%%%%%%%%%%%%%%%%%%%%%%%%%%%%%%%

\section{Proofs of Theorem \ref{thm:counter} and Corollaries \ref{cor:torus} 
and \ref{cor:all_unk_nums}\label{sec:counter}}

%%%%%%%%%%%%%%%%%%%%%%%%%%%%%%%%%%%%%%%%%%%%%%%%%%%%%%%%%%%%%%%%%%%%%%%%%%%%
%%%%%%%%%%%%%%%%%%%%%%%%%%%%%%%%%%%%%%%%%%%%%%%%%%%%%%%%%%%%%%%%%%%%%%%%%%%%

The two summands which make up our example in Theorem \ref{thm:counter} are $K=7_1$ and
its mirror image $\overline{K}$. $K$ is the $(2,7)$-torus knot. A routine 
calculation establishes that $K$ has signature
$\sigma(K)=6$; from the inequality $|\sigma(K)|/2\leq u(K)$ \cite{mu65}
we know that $K$ has unknotting number at least 3. 
Since any 3 crossing changes in its 
standard 7-crossing diagram results in the unknot, we have $u(K)=3$. Alternatively, 
one can appeal to the important
and more far-reaching result of Kronheimer and Mrowka \cite{krmr} that
the unknotting number of the $(p,q)$-torus knot $T(p,q)$ is $(p-1)(q-1)/2$.
It is a straightforward result that unknotting number is preserved under mirror image
- hold a mirror up to the unknotting sequences - and so $u(\overline{K})=3$,
as well. Consequently, $u(K)+u(\overline{K})=3+3=6$.

We will demonstrate, however, that $L=K\#\overline{K}$ has unknotting number at most 
$5$. We do this by giving a diagram of $L$ for which changing
two crossings gives (a diagram for) the knot $K_1=K14a18636$, 
in the notation from SnapPy. We then show that the knot $K_2=K15n81556$
has unknotting number (at most) 2, and has a diagram for which one
crossing change yields a diagram for $K_1$. 
Consequently, $L=K\#\overline{K}$ has a diagram for which 2 crossing changes
can take us to a diagram of $K_1$; an isotopy then takes us to 
another diagram for $K_1$ that admits a crossing change to a diagram
of $K_2$, which has unknotting number at most $2$.
So it takes at most $2+1+2=5$ crossing changes, with
isotopies, to transform our diagram of $K\#\overline{K}$ to the unknot.
This unknotting sequence therefore has length $5$ (that is, $6$ knots and $5$ 
crossing changes);
consequently, $u(K\#\overline{K})\leq 5$.

\medskip

\begin{figure}[h]
\begin{center}
\includegraphics[width=5in]{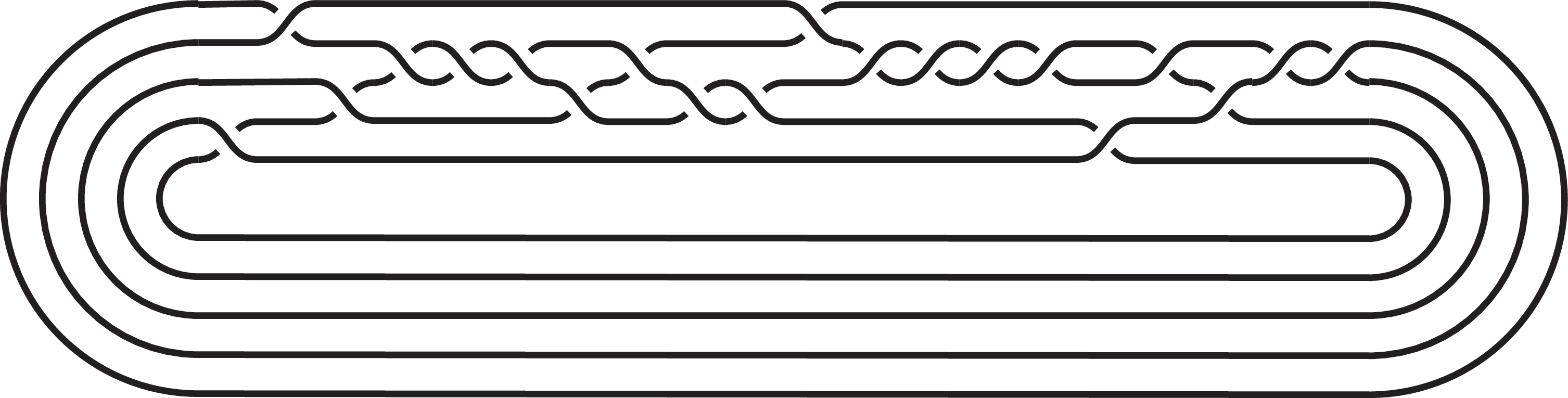}
\caption{The initial braid closure}\label{fig:braid}
\end{center}
\end{figure}

The initial diagram for $L$ in our unknotting sequence was found as the closure of a braid,
using the program SnapPy \cite{snappy}.
For our example, the particular braid ultimately produced, which yields the diagram for the first 
half of our unknotting sequence, is given by the braid word:

\begin{center}
$[\underbar{1}, \underbar{-4}, 2, 3, 3, 3, 2, 3, 2, 2, 4, -3, -3, -3, -3, -1, -3, -2, -3, -3]$
\end{center}

\noindent which represents a knot diagram with 20 crossings. Here we use
SnapPy's convention that standard braid generators are listed by their
(signed) subscripts; a positive entry means a generator, negative means its inverse.
This braid closure is given in Figure \ref{fig:braid}; the braid is read left to right, 
with braid generator $1$ at bottom and $4$ at top.

\begin{figure}[h]
\begin{center}
\includegraphics[width=5in]{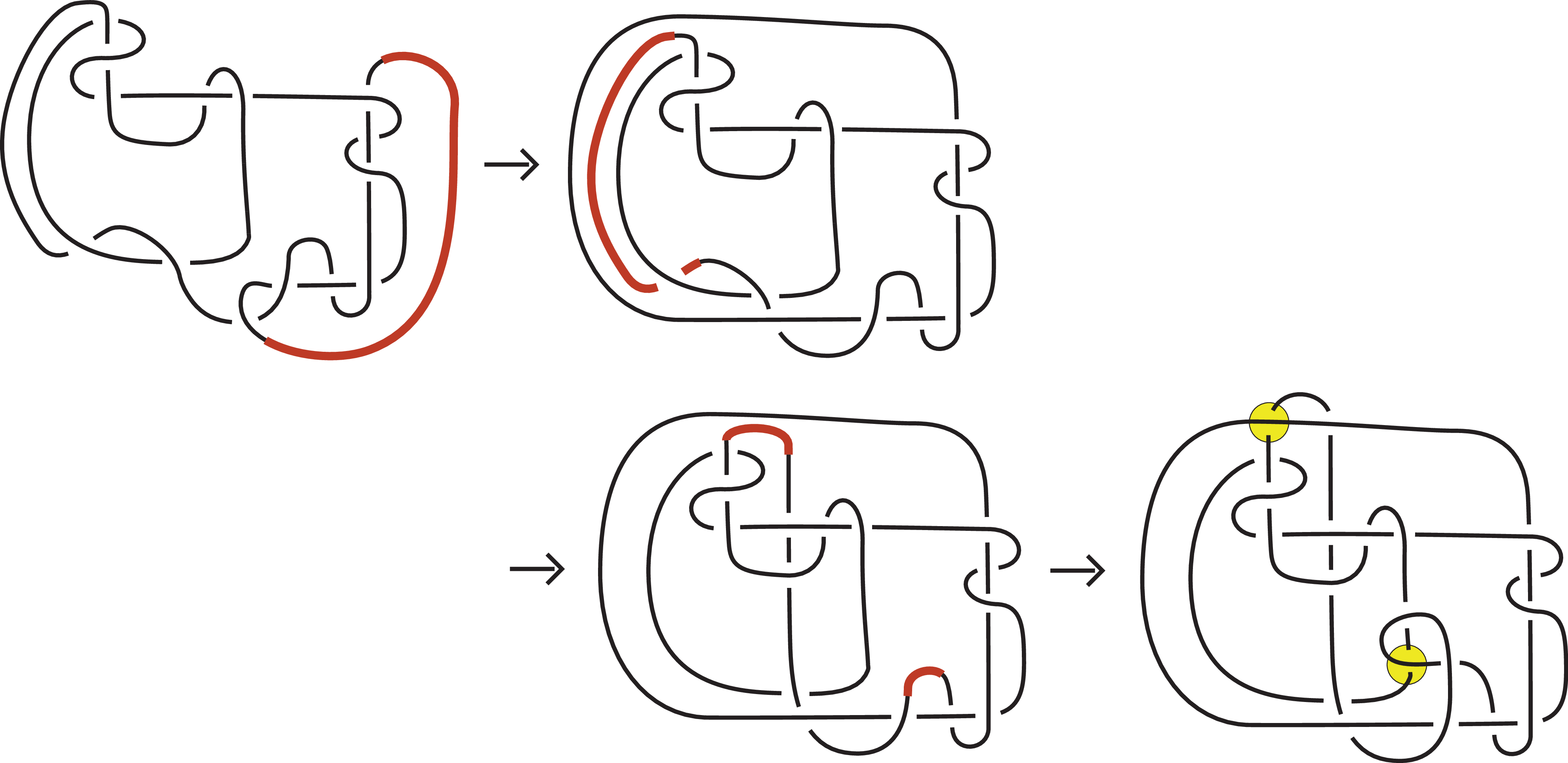}
\caption{Diagrams for $7_1\#\overline{7_1}$}\label{fig:7171}
\end{center}
\end{figure}

The reader can verify in SnapPy/Sage \cite{snappy},\cite{sage}
that this is a braid representative 
for $L$ by using the command {\tt K=Link(braid\_closure=brd)}
or {\tt K=snappy.Link(braid\_closure=brd)} (calling the above tuple `brd')
and employing the {\tt simplify()} and {\tt deconnect\_sum()} commands
to recover the summands $7_1$ and $\overline{7_1}$ . Alternatively, one can 
start from a diagram for 
$L=7_1\#\overline{7_1}$, on the left side of Figure \ref{fig:7171}, and arrive at 
the diagram of the knot that SnapPy will produce from the 
braid word above, after a few Reidemeister moves, as indicated in Figure \ref{fig:7171}.
Strands that are not fixed in each step of this isotopy, and in the isotopies that
are represented in the diagrams 
that follow, are represented as thicker (and red) arcs.

Crossing changes on this diagram can be carried out by changing the signs
of entries in the tuple, to create a new braid word representing the
knot after crossing change. In the above braid word we have underlined the entries
relevant to our construction. Following the convention that SnapPy/python starts labeling entries
in a tuple at `0', changing the crossings at indices 0 and 1 yields the 
braid word 

\begin{figure}[h]
\begin{center}
\includegraphics[width=5in]{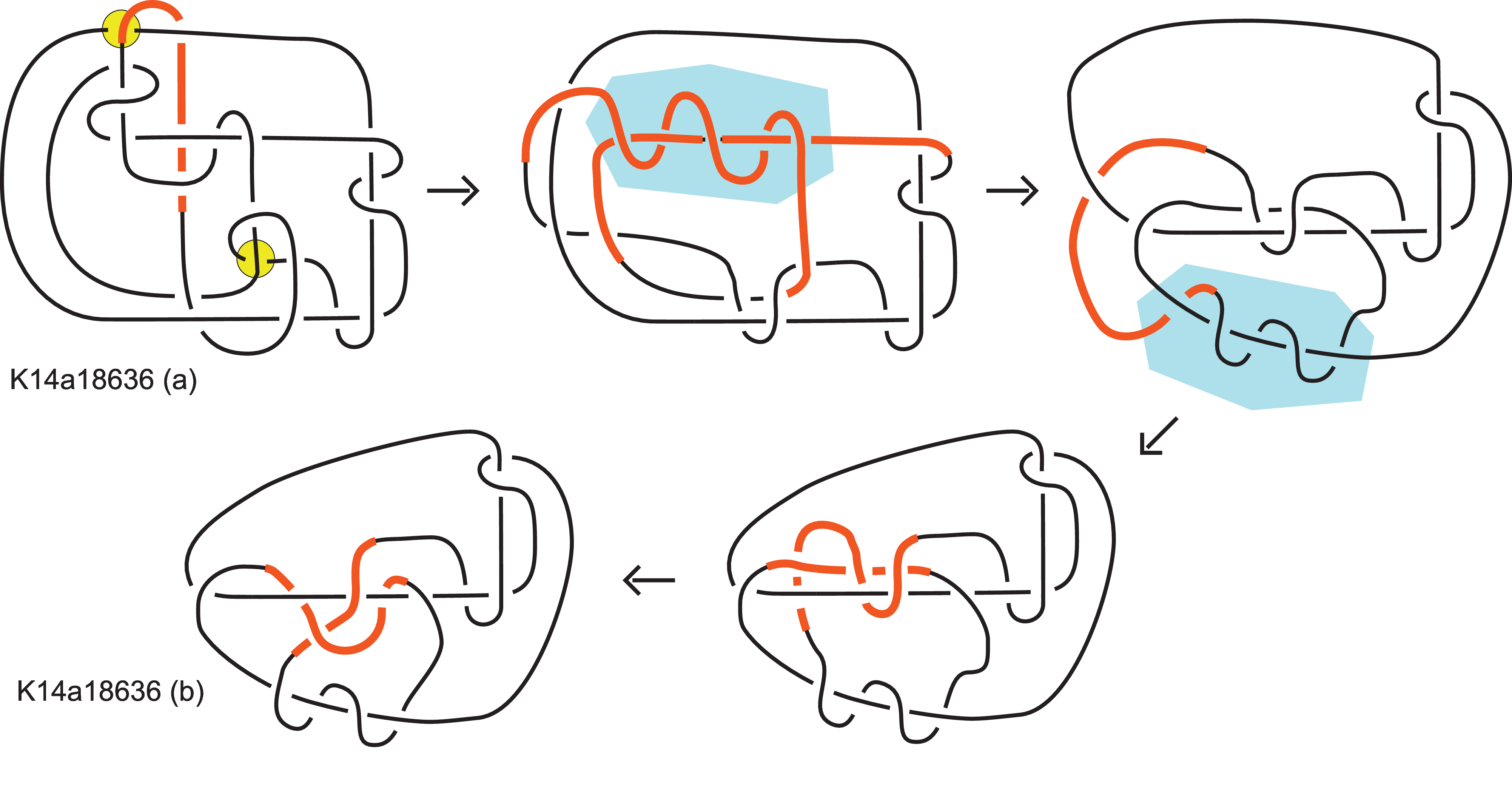}
\caption{Diagrams for $K14a18636$}\label{fig:K14a18636}
\end{center}
\end{figure}

\begin{center}
$[\underbar{-1}, \underbar{4}, 2, 3, 3, 3, 2, 3, 2, 2, 4, -3, -3, -3, -3, -1, -3, -2, -3, -3]$
\end{center}

\noindent which SnapPy identifies as belonging to the knot $K_1=K14a18636$. 
These crossing changes, and subsequent crossing changes in further diagrams, 
are marked with a shaded (and yellow) circle in Figure \ref{fig:7171}. A diagram 
for this knot $K_1$ built from 
the closure of the braid represented by this braid word 
is given in Figure \ref{fig:K14a18636};
the two crossing changes made to obtain it from the rightmost diagram in 
Figure \ref{fig:7171} are circled.
We also provide a schematic
of an isotopy from this 20-crossing projection to a 14-crossing projection of $K_1$.
Since unknotting number changes by at most one under crossing change, 
and we can recover $L$ from $K_1$ by reversing these $2$ crossing changes, we have
$u(L)\leq u(K_1)+2$ .

We now demonstrate that $K_1$ has unknotting number at most $3$.
We do this by making use of a pair of 15-crossing diagrams of the knot $K_2=K15n81556$.
For this discussion we utilize Dowker-Thistlethwaite (DT) codes \cite{dt83} for the knots, 
where, again, each entry of the code represents a crossing in the diagram, 
and a crossing change can be achieved by changing the sign of an entry of the code.
The knot $K_2$ has, from the authors' earlier work on the Bernhard-Jablan conjecture
\cite{bh21}, 100 distinct 15-crossing projections. For one of them, with DT code

\begin{center}
$[\underbar{-4},-16,24,26,18,20,28,22,-2,10,12,30,6,8,14]$,
\end{center}

\noindent a single crossing change, at index 0, yields the DT code

\begin{figure}[h]
\begin{center}
\includegraphics[width=5in]{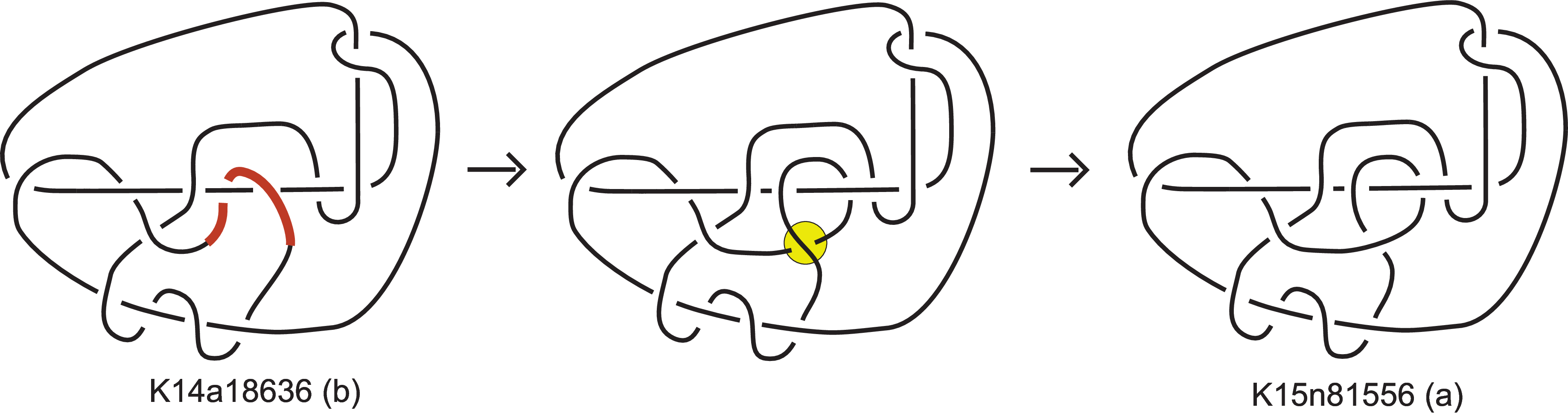}
\caption{$K14a18636$ to $K15n81556$ (a)}\label{fig:K1_to_K2}
\end{center}
\end{figure}

\begin{center}
$[\underbar{4},-16,24,26,18,20,28,22,-2,10,12,30,6,8,14]$,
\end{center}

\noindent which SnapPy identifies as belonging to the knot $K_1=K14a18636$. 
Figure \ref{fig:K1_to_K2} shows how the last diagram in Figure \ref{fig:K14a18636},
for $K14a18636$, can be transformed via Reidemeister moves and a single crossing change
into the diagram for $K15n81556$ corresponding to this DT code.

On the other hand, a different diagram for $K_2$ is given by the DT code 

\begin{center}
$[4,12,-24,14,18,2,\underbar{20},26,8,10,-28,-30,16,-6,-22]$, 
\end{center}

\begin{figure}[h]
\begin{center}
\includegraphics[width=5in]{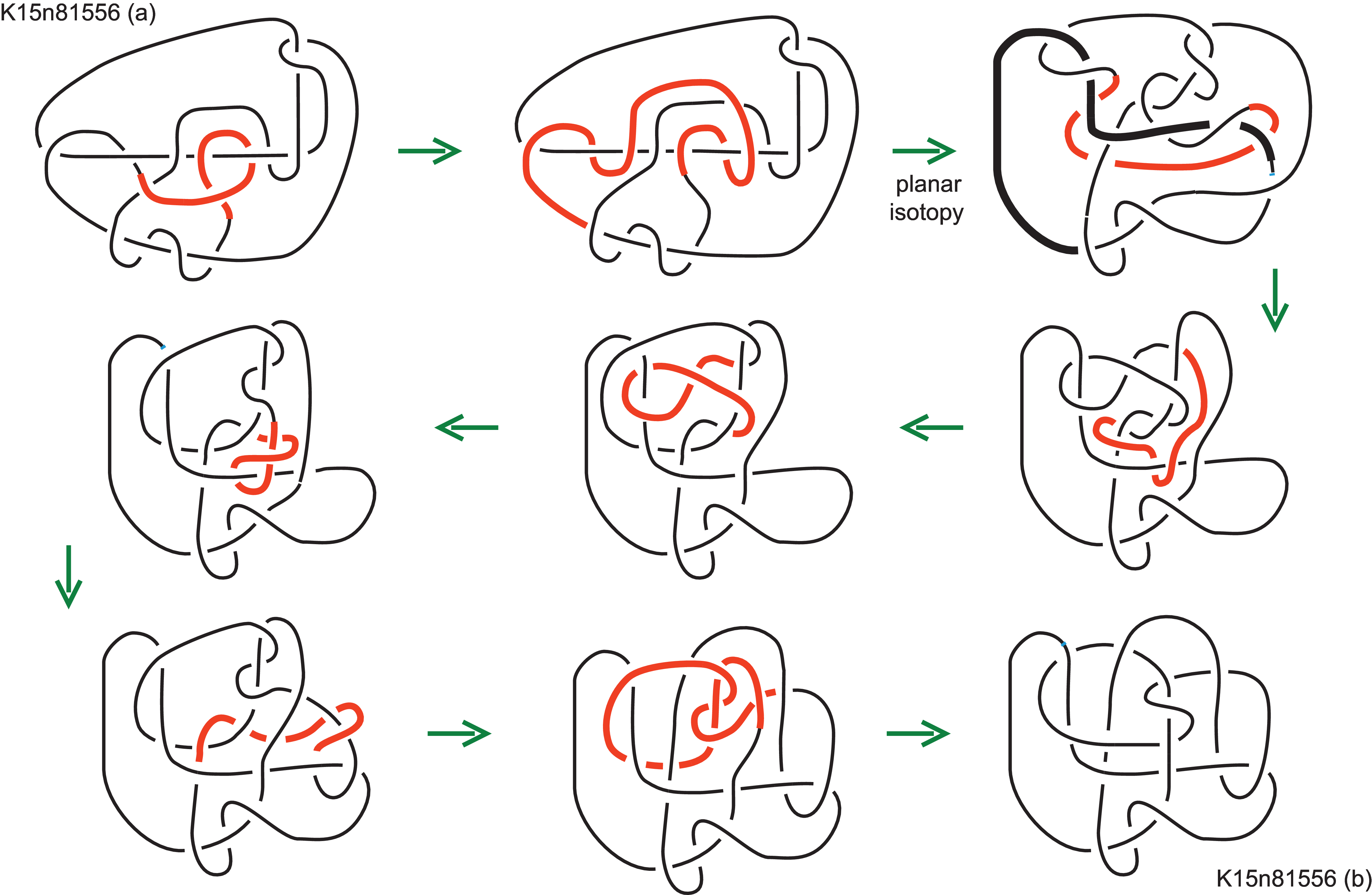}
\caption{$K15n81556$ (a) to $K15n81556$ (b)}\label{fig:K15n_to_K15n}
\end{center}
\end{figure}

\noindent A sequence of isotopies interpolating between the two diagrams
for $K_2$ is shown in 
Figure \ref{fig:K15n_to_K15n}; a different such sequence can be found in \cite{wz25}.
In the figure, `planar isotopy' means that there are no changes in the 
combinatorics of the underlying 
knot diagram; the highlighted arc is straightened to the corresponding arc marked.

\noindent When we change a single crossing, at index 6, this yields the DT code

\begin{center}
$[4,12,-24,14,18,2,\underbar{-20},26,8,10,-28,-30,16,-6,-22]$, 
\end{center}

\noindent which SnapPy identifies as belonging to the knot $K_3=K12n412$. The knot 
$K_3$ has unknotting 
number 1 (as noted, for example, on the Knotinfo site \cite{knotinfo}); 
and in fact changing the crossing in the DT code above for $K_3$, at index 13, yields the code

\begin{figure}[h]
\begin{center}
\includegraphics[width=5in]{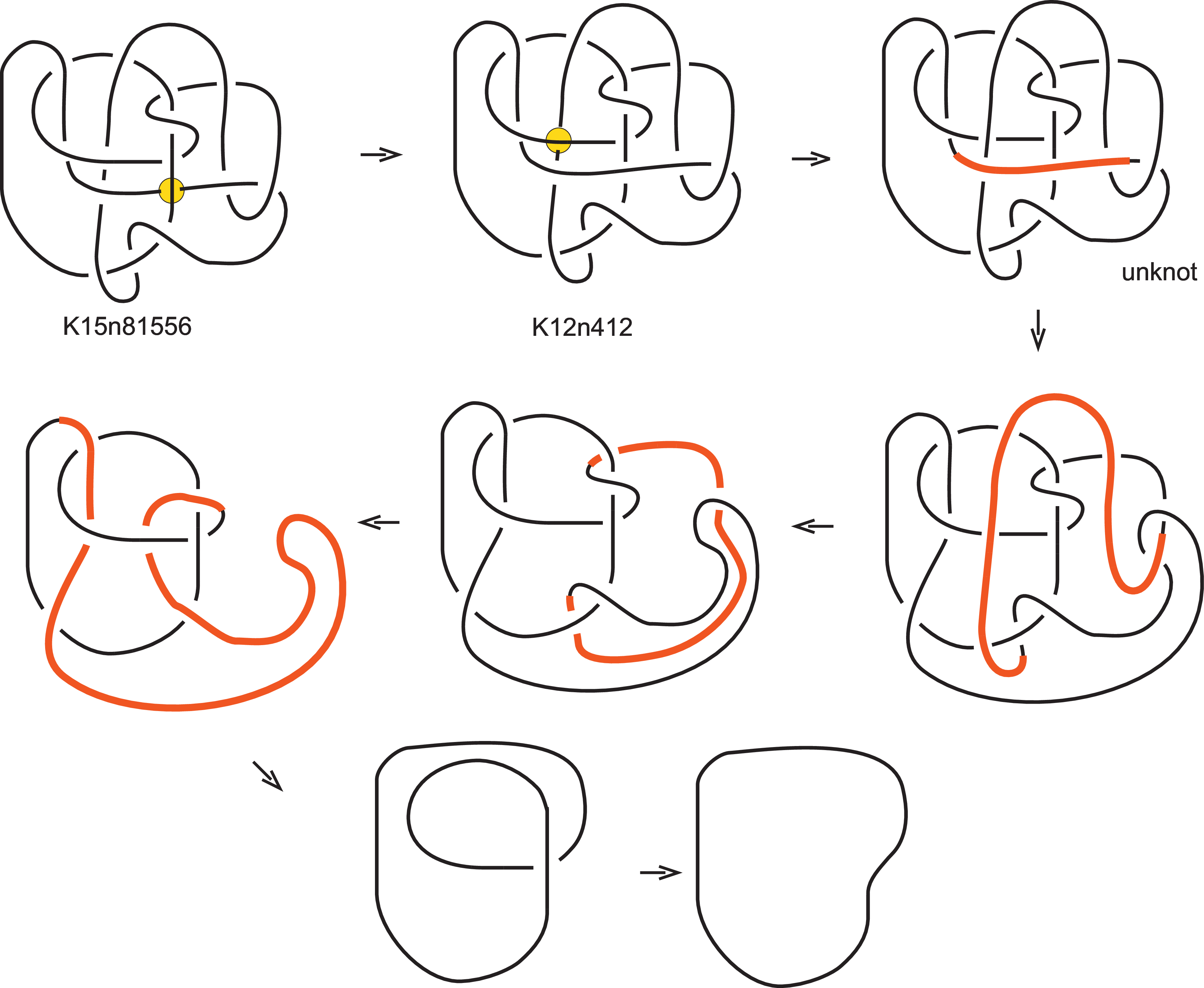}
\caption{$K15n81556$ (b) to the unknot}\label{fig:K15_to_unknot}
\end{center}
\end{figure}

\begin{center}
$[4,12,-24,14,18,2,-20,26,8,10,-28,-30,16,\underbar{6},-22]$,
\end{center}

\noindent which is a diagram of the unknot. 
A sequence of isotopies from the last diagram for $K15n81556$ in Figure \ref{fig:K15n_to_K15n}
to the unknot, after making these two crossing changes, 
is shown in Figure \ref{fig:K15_to_unknot}. Consequently, we can obtain the 
unknot starting from $K_1=K14a18636$ with three crossing changes; 
first to $K_2=K15n81556$, then to $K_3=K12n412$, and finally to the unknot.

\medskip

This establishes that $u(K_1)=u(K14a18636)\leq 3$. 
Therefore, 

\begin{center}
$u(7_1\#\overline{7_1})=u(L)\leq u(K_1)+2\leq 3+2=5$ .
\end{center}

\noindent Since $u(7_1)+u(\overline{7_1})=3+3=6$, this shows
that $u(7_1\#\overline{7_1}) < u(7_1)+u(\overline{7_1})$, establishing the 
failure of additivity for these knots. This completes the proof 
of Theorem \ref{thm:counter}. \hfill $\square$

\bigskip

It is a well-known fact that in computing unknotting number,
crossing changes and isotopies can, in a certain sense, be exchanged,
at the expense of the crossing changes taking place in a more
complicated diagram. That is, every knot $K$ has a diagram $D$ so that
$u(K)$ crossing changes in that diagram yields a diagram for the unknot.
After the authors posted an earlier version of this paper to 
arxiv.org, using some of the sequences of crossing changes and 
isotopies that we found above 
Wang and Zhang \cite{wz25} found
such a diagram $D$ for $7_1\#\overline{7_1}$ for which $5$ crossing 
changes in the diagram yields a diagram of the unknot. Their 
diagram for $7_1\#\overline{7_1}$ has $56$ crossings.

\smallskip

Next we turn our attention to the proof of Corollary \ref{cor:torus},
extending our main example to an infinite class of examples.
A knot $K_1$ is called \emph{Gordian adjacent} to the knot $K_2$ \cite{fel14} if 
$u(K_2)>u(K_1)$ and there are 
$u(K_2)-u(K_1)$ crossing changes on a diagram for $K_2$ that results 
in a diagram for $K_1$. Put differently, $K_2$ possesses a minimal unknotting 
sequence which contains the knot $K_1$.
For any pair of knots $L_1,L_2$ that the knot $7_1$ is Gordian adjacent to,
we can build an unknotting sequence for $L_1\#\overline{L_2}$ using 
fewer than $u(L_1)+u(L_2)$ crossing changes. In particular, starting with $L_1\#\overline{L_2}$,
applying $u(L_1)-3=u(L_1)-u(7_1)$ crossing changes to the $L_1$-summand, and 
$u(L_2)-3$ crossing changes to the $\overline{L_2}$ summand will yield a diagram for
$7_1\#\overline{7_1}$, and this knot can be unknotted with $5$ crossing changes (after isotopy). 
Consequently, $(u(L_1)-3)+(u(L_2)-3)+5=u(L_1)+u(L_2)-1$ crossing changes suffice to unknot
$L_1\#\overline{L_2}$, and so $u(L_1\#\overline{L_2})\leq u(L_1)+u(L_2)-1<u(L_1)+u(L_2)$.

One family of knots which can quickly be seen to satisfy this condition is the 
set of torus knots $T(2,2k+1)$ with 
$k\geq 3$. These knots have unknotting number $(2-1)((2k+1)-1)/2=k$, and changing
any $k-3$ crossings in their standard diagrams will yield $T(2,7)=7_1$
(since changing a crossing in the standard diagram for $T(2,N)$ will, after a 
Reidemeister II move, yield the standard diagram for $T(2,N-2)$, and we are then 
done by induction). This yields Corollary \ref{cor:torus}. \hfill $\square$

\bigskip

Further extensions of Theorem \ref{thm:counter} can be made to pbtain
more examples in a similar way. 
Several authors have explored the construction of minimal unknotting sequences for torus knots
\cite{ba10},\cite{sima15}, and the problem of Gordian adjacency for torus 
knots \cite{fel14},\cite{bl16}. Feller's work \cite{fel14} shows that 
all torus knots $T(p,q)$ except (possibly) $T(2,3),T(2,5), T(3,4), T(4,5)$ and $T(5,6)$
have $7_1=T(2,7)$ Gordian adjacent to them. Tetsuya Abe \cite{abe25} has pointed
out that work of Kanenobu \cite{kan10} shows that $T(2,7)$ is Gordian adjacent to 
both $T(4,5)$ and $T(5,6)$, and has provided the authors with explicit unknotting
sequences to show this. $T(2,7)$ cannot be Gordian adjacent to the 
knots $T(2,3),T(2,5)$ and $T(3,4)$, however, because their unknotting numbers are too low; 
they have unknotting numbers 1, 2, and 3, respectively. Taken together,
these results show that nearly every torus knot can be used to build a
counterexample to the additivity conjecture.

\begin{corollary}\label{cor:most_torus}
For the collection of torus knots

${\mathcal T}=\{T(p,q)$\ : $2\leq p<q,\ \gcd(p,q)=1,\ (p,q)\not\in\{(2,3),(2,5),(3,4)\}\}$, 

\noindent any pair of knots $K_1,K_2\in{\mathcal T}$ satisfy
$u(K_1\#\overline{K_2})<u(K_1)+u(K_2)$.
\end{corollary}

It is an interesting open question to decide if the torus knots $T(2,3)$, $T(2,5)$ or $T(3,4)$
can be paired with other knots to form a connected sum which also fails additivity.

For an example that is not a torus knot, we can use the knot $10_{139}$ (which, 
in fact, was a summand of the actual first example that the authors found in 
proving Theorem \ref{thm:counter}; see Section \ref{sec:howfound} below).  
The knot $L_1=10_{139}$ is a non-alternating hyperbolic knot; it was shown by Kawamura \cite{kawa}
that $u(L_1)=4$. In fact, Kawamura showed, even more, that $L_1$
has smooth 4-genus $g_4(L_1)=4$ (and $g_4(L)\leq u(L)$ for every knot $L$ \cite{mu65}); 
the 4-genus of $L_1$ was also established by different means (using their tau invariant) 
by Ozsv\'ath and Szab\'o \cite{os03}. One DT code for $10_{139}$ is given by

\begin{center}
[12, 14, -10, -20, \underbar{-16}, 18, 2, -8, 4, -6],
\end{center}

\noindent and changing the sign of the entry at index 4 yields the DT code

\begin{center}
[12, 14, -10, -20, \underbar{16}, 18, 2, -8, 4, -6],
\end{center}

\noindent which belongs to the knot $7_1$. As a result, an unknotting sequence for $7_1$
(with $3$ crossing changes) can be augmented with this crossing change to yield
an unknotting sequence for $10_{139}$ with $4$ crossing changes, that is, a
minimal unknotting sequence. Consequently, $T(2,7)$ is Gordian adjacent to
$10_{139}$, yielding:

\begin{corollary} $u(7_1\#\overline{10_{139}})<u(7_1)+u(10_{139})$ and 
$u(10_{139}\#\overline{10_{139}})<u(10_{139})+u(10_{139})$ .
\end{corollary}

We expect that there are many other knots whose unknotting numbers are
known and that contain $7_1$ in
a minimal unknotting sequence, which would then yield other connected
sums that fail additivity. In particular, any intermediate knot $K$
in a Gordian adjacency from the torus knot $T(p,q)$ to the knot $7_1$ also
has $7_1$ Gordian adjacent to it, and so can be used as a summand in such examples.

For further examples, by leveraging a result of Baader \cite{ba06}, we can build 
a great many non-torus
knots for which $u(L_1\# \overline{L_2})<u(L_1)+u(L_2)$. Baader shows that for any two knots 
$K$, $L$ that are two crossing
changes apart, there are infinitely many distinct knots $K_i$ that are intermediate
to them, in the sense that they are each one crossing change away from 
both $K$ and $L$. By finding infinitely
many such knots $K_i$ intermediate between the knots 
$T(2,2(n+1)+1)$ and $T(2,2(n-1)+1$, for $n\geq 4$, we can be certain of two things:
First, each of the knots $K_i$ has unknotting number precisely $n$, 
since they are each one crossing change away from 
$T(2,2(n+1)+1)$, with unknotting number $n+1$, and from $T(2,2(n-1)+1)$, 
with unknotting number $n-1$. The value of $u(K_i)$
must then be at most one away from both $n+1$ and $n-1$. Second, $7_1$ is Gordian
adjacent to each $K_i$, since, in our standard unknotting sequence for
$T(2,2(n+1)+1)$, passing through $T(2,2(n-1)+1$ and $7_1$, we can replace
the knot $T(2,2n+1)$ with $K_i$ and still have an unknotting sequence
of the same length, which passes through $K_i$ and then $7_1$. 
This gives us Corollary \ref{cor:all_unk_nums}.  \hfill $\square$

\bigskip

Every diagram $D$ of a knot $K$ can be turned into a diagram of the unknot by changing 
at most $\lfloor c(D)/2\rfloor$ crossings, where $c(D)$ is the number of 
crossings in the diagram. This can be done, for example, by changing crossings to build an 
ascending diagram of the unknot from the diagram $D$. This leads to Ozawa's result \cite{oza10}
that $u(K)\leq\lfloor (c(K)-1)/2\rfloor$,
where $c(K)$ is the \emph{crossing number} of $K$, the minimum number of crossings
among all diagrams for $K$. It was shown by Taniyama \cite{tan09} that only 
the torus knots $K=T(2,n)$ satisfy $u(K)=\lfloor (c(K)-1)/2\rfloor$.
Abe, Hanaki, and Higa \cite{ahh12} showed that if $K$ is prime, 
then $u(K)=\lfloor (c(K)-2)/2\rfloor$
if and only if $K$ is either the knot $4_1$ (the figure-8 knot)
or $K$ can be represented as a positive or negative 3-braid. 
They also showed that the only possible
non-prime knots with $u(K)=\lfloor (c(K)-2)/2\rfloor$ are the knots $T(2,n)\# T(2,m)$ and
$T(2,n)\# \overline{T(2,m)}$ for $n,m\geq 3$, and that the 
knots $K=T(2,n)\# T(2,m)$ do satisfy this equality. They
conjectured that the knots $T(2,n)\# \overline{T(2,m)}$ also satisfy 
$u(K)=\lfloor (c(K)-2)/2\rfloor$, and showed that this is equivalent to 
the statement that $u(T(2,n)\# \overline{T(2,m)})=u(T(2,n))+u(\overline{T(2,m)})$.
Our Corollary \ref{cor:torus} shows that this last assertion is false, when $n,m\geq 7$.

This leaves the cases $T(2,3)\# \overline{T(2,m)}$ 
and $T(2,5)\# \overline{T(2,m)}$ for $m\geq 3$ and odd, 
to decide whether or not $u(K)=\lfloor (c(K)-2)/2\rfloor$
(that is, whether or not 
$u(T(2,3)\# \overline{T(2,m)})=1+\frac{m-1}{2}$ 
and
$u(T(2,5)\# \overline{T(2,m)})=2+\frac{m-1}{2}$).
Settling these last two cases would complete the characterization of the knots $K$ with
$u(K)=\lfloor (c(K)-2)/2\rfloor$ .

%%%%%%%%%%%%%%%%%%%%%%%%%%%%%%%%%%%%%%%%%%%%%%%%%%%%%%%%%%%%%%%%%%%%%%%%%%%%
%%%%%%%%%%%%%%%%%%%%%%%%%%%%%%%%%%%%%%%%%%%%%%%%%%%%%%%%%%%%%%%%%%%%%%%%%%%%

\section{How (and why) the main example was found}\label{sec:howfound}

%%%%%%%%%%%%%%%%%%%%%%%%%%%%%%%%%%%%%%%%%%%%%%%%%%%%%%%%%%%%%%%%%%%%%%%%%%%%
%%%%%%%%%%%%%%%%%%%%%%%%%%%%%%%%%%%%%%%%%%%%%%%%%%%%%%%%%%%%%%%%%%%%%%%%%%%%

The process of finding the diagram for $7_1\#\overline{7_1}$ above and the discovery
of its unexpectedly low unknotting number came as a consequence of another 
ongoing project of the authors: our search
for an example to show that unknotting number is not
invariant under mutation. This would answer
Problem 1.69(C) in Kirby's problem list \cite{kir93}, in the negative. 
Our approach to that problem included
building a collection of knots which would have unknotting number 4 (if unknotting 
number were additive), 
but which would have comparatively small signature. 
The inequality $|\sigma(K)|/2\leq u(K)$ \cite{mu65} is useful in bounding unknotting 
number from below, but it can be difficult to find knots
where the difference between these two quantities is known to be greater than 1
(which we needed for our search), although work of Stoimenow \cite{sto02} is a significant
outlier, which deserves to be more widely known.
On the other hand, the facts that signature is additive under connected sum,
and that the signature of the mirror image of a knot $K$ is the negative of the
signature of $K$ (see, e.g., \cite{mu65}), means that is it a routine matter to 
build connected sums where
the unknotting number would be large, if unknotting 
number were additive, but where the signature 
is small. It was for this reason that we began to search for knots a few
crossing changes away from connected sums of low crossing number knots with 
known unknotting number, in order to build a collection of new knots
with (if unknotting number were additive) high enough unknotting number
to be useful for our project.

Our approach to finding the new knots was to adopt the process we had 
followed for two previous projects, 
where we settled the Bernhard-Jablan Conjecture (in the negative) \cite{bh21}
and where we completed the computation of the smooth 4-genus of the prime knots
through 12 crossings \cite{bh22}. We started with pairs 
of knots, one with known unknotting number 3 and the other 
with known unknotting number 4, took their connected sum, and then built
diagrams for these sums for which crossing changes stood a good chance of
building prime knots, with low enough (after simplification) crossing number that SnapPy
could recognize them (that is, crossing number 15 or less). We chose unknotting 
numbers 3 and 4 so that we were starting with connected sums that were likely to have
unknotting number 7, and then randomly changed three crossings so 
that the resulting knots conjecturally would have unknotting number (at least) 4.
We built the knot diagrams by starting with the standard diagrams
for the connected sums, randomly
performing some crossing-number-increasing Reidemeister moves
using SnapPy's {\tt backtrack()} command, and
then converting them to braids using
the {\tt braid\_word()} command. For this last step, SnapPy uses Vogel's algorithm \cite{vog90}
to convert the knot diagram into a closed braid. We chose this approach 
because (a) SnapPy can quickly carry
out both of these operations and (b) our experience with the previous two 
projects had shown us that working with braid representatives
built in this way succeeds in building diagrams that 
are randomized enough to yield good results (that is, many different knot types)
when we change crossings in the braid representations. 
Experience also showed that changing 3 crossings, when starting from a 
connected sum in this way, resulted in 
a good proportion of the knots that were generated being prime; 2 crossing changes was typically 
insufficient. Using braid words, changing crossings is very quick; it amounts 
to changing the signs of entries in the braid representation.

Our work in \cite{bh21} had provided
upper bounds on the unknotting number of all of the prime knots through
15 crossings, what we called, in that paper, their \emph{weak BJ-unknotting number},
that is, their unknotting number \underbar{if} the Bernhard-Jablan conjecture 
were true.
We could have confidence that the vast majority of
these computations represented the actual unknotting number of the knots, 
based on comparison with known values for lower crossing number knots.
We compared the knots coming from the randomized procedure described above to our list of
`weak BJ-unknotting number 4' knots to create a group of what we called 
`conditionally unknotting number 
4' knots. Working with several million braids yielded 240 such 
14- and 15-crossing knots. 
One of these was the knot $K14a18636$.

We then began to reverse the search, using these 240 knots as 
initial knots to carry out crossing changes on their braid representatives,
in the same manner in which the connected sum diagrams had been built.
For the mutation project, our goal was to 
show that one of a collection of target knots had unknotting number 3, 
by showing that it was crossing-adjacent to a knot with unknotting number 4.
These target knots came from a database of mutant knot pairs constructed by
by Stoimenow \cite{stoi}, using our Bernhard-Jablan data to identify
mutant pairs with distinct weak BJ-unknotting number.
Prior such searches using knots with unknotting number known to 
be 4 had failed to make such a connection, and so our plan was to search 
among knots one crossing change away from these 240 knots with conditional 
unknotting number 4, instead. Connecting one of the 240 initial knots 
with one of these mutant
`targets' would have proved that either unknotting number
is not invariant under mutation, or that the knot from among our 240
that made the connection would actually have unknotting number at most 3, 
and so, as argued in Section \ref{sec:counter}, unknotting number 
would not be additive under connected sum.
(But, were such a thing to happen, we would not know which of the two possibilities was true!)

At the start of this phase of the project it was not clear 
what this process would produce. The most likely prospect seemed to be that
we would never find a connection between our 240 initial knots and our list of 
target knots, if only because both conjectures might be true. Our 
expectation, if the project succeeded, was that it would be because 
the target knot we connected with had unknotting number 3, and so
unknotting number would not be mutation invariant. It came as quite a
surprise then when, instead, this reverse search
found a crossing change from one of the 240 knots that we had identified, 
namely $K14a18636$, to the knot $K15n81556$, which had a known, 
and unexpectedly low, weak BJ-unknotting number of 2. 

The first example we found started from a closed braid representation of
the knot $7_1\#10_{139}$;
three crossing changes, that our search had found, transformed this 
to a representative of the 
knot $K14a18636$.
One of those three crossing changes transformed the braid
into a braid representative for the knot 
$7_1\#\overline{7_1}$, and so this braid word contained two 
crossing changes to a representative of the knot $K14a18636$. Initially, 
the braid word that made this connection had length 119, and so 
represented a 119-crossing knot diagram. 
But several simplifications of the
braid closure could, it turned out, be carried out. 
These mostly involved Reidmeister I moves `hidden' by a sequence
of Reidemeister II moves. Their discovery was aided by the program KnotPlot
\cite{knotplot}, to 
visualize long stretches of the braid and find the crossings to eliminate. 
This enabled the authors to greatly
reduce the number of crossings needed (from 119 to 20), and the number of 
strands in the braid needed (from 14 to 5). This resulted in the diagram
and braid word used in our discussion in Section \ref{sec:counter}.

This example, of course, 
destroyed the hypothesis underlying our mutation project;
but that would seem to be a small price to pay, for the result that it achieved.

%%%%%%%%%%%%%%%%%%%%%%%%%%%%%%%%%%%%%%%%%%%%%%%%%%%%%%%%%%%%%%%%%%%%%%%%%%%%
%%%%%%%%%%%%%%%%%%%%%%%%%%%%%%%%%%%%%%%%%%%%%%%%%%%%%%%%%%%%%%%%%%%%%%%%%%%%

\section{The road forward}\label{sec:road}

%%%%%%%%%%%%%%%%%%%%%%%%%%%%%%%%%%%%%%%%%%%%%%%%%%%%%%%%%%%%%%%%%%%%%%%%%%%%
%%%%%%%%%%%%%%%%%%%%%%%%%%%%%%%%%%%%%%%%%%%%%%%%%%%%%%%%%%%%%%%%%%%%%%%%%%%%

With the additivity of unknotting number being settled in the negative, 
the quest to build a more comprehensive understanding 
of unknotting number becomes more challenging. The fact that
unknotting number does not behave as expected with respect to connected sum
means that building a purely computational approach to determining
unknotting numbers, as one could argue the other invariants mentioned in the
introduction enjoy, seems further out of reach. In particular, we cannot 
rely on knowledge of unknotting numbers for prime knots to provide all of
the information we need to compute the invariant more generally. 

As noted in the proof of Corollary \ref{cor:knotgroup} in Section \ref{sec:intro}, 
since the knot group of a connected sum is insensitive to 
mirroring, but its unknotting number can change depending on the mirroring
of the summands, this means that the knot group alone is 
not sufficient to detect the unknotting number of knots. 
On the other hand, the knot group $\pi(K)=\pi_1(S^3\setminus K)$ 
together with the knot's peripheral structure
(the conjugacy class of the image of $\pi_1(\partial N(K))$ in $\pi(K)$ 
under the inclusion-induced homomorphism)
does determine the knot $K$ \cite[Corollary~6.5]{wa68}, and therefore (in the abstract) 
determines its unknotting number. Any general method to compute the unknotting number
of a knot from data coming from the knot group would therefore need to 
take into account this peripheral structure. One approach might be to focus on
the 2-generalized knot group $G_2(K)$, which is built from a presentation of the
knot group and which is known \cite{nn08} to determine the knot, up to mirroring.
It therefore, in the abstract, determines the unknotting number, as well.

\begin{question}
How does $G_2(K)$ behave under connected sum? Can $u(K)$ be computed from $G_2(K)$?
\end{question}

Knowing that additivity of unknotting number can in
fact fail makes it more interesting to discover when it does not
fail. That is, establishing that $u(K\#K^\prime)=u(K)+u(K^\prime)$
for families of knots has increased relevance and importance, since we can
no longer expect this to hold in general. As mentioned in 
Section \ref{sec:intro}, there are families
${\mathcal K}$ of knots for which additivity holds for pairs of knots
chosen within the family \cite{krmr},\cite{ba05}; further classes would 
help give us a better picture of where we might look for more examples
of non-additivity. This in turn might lead to new shortcuts, which 
would have relevance to the computation of unknotting number more generally.

Since additivity of unknotting number does not hold in general, 
we can also explore how badly additivity can fail. 
We can leverage the examples above to show that 
the difference $u(K\# L)-u(K)-u(L)$ can be arbitrarily large, by, for example, using
$K=\#_n 7_1$ and $L=\#_n \overline{7_1}$, the $n$-fold connected
sums of these knots. We have $u(K)=u(L)=3n$
by signature considerations, but by rearranging terms we have
$K\#L=\#_n(7_1\#\overline{7_1})$, with unknotting number at most $5n$. 
Any knot that the knot $7_1$ is gordian adjacent to can be used in
same fashion, and so the knots described in Corollaries 
\ref{cor:all_unk_nums} and \ref{cor:most_torus} can be used to build infinitely
many examples. This gives Corollary \ref{cor:large_gaps}.   \hfill $\square$

\bigskip

It is an interesting question to 
ask to what extent such behavior also occurs for prime knots $K$ and $L$. 
Sebastian Baader \cite{ba25} has pointed out that the pretzel knot $L=P(7,7,\ldots,7,2)$
with $n$ 7's,
which is prime, can be transformed to the knot $\#_n 7_1$ by a single crossing
change in the last collection of half-twists. This implies that $u(L\#\overline{L})\leq 5n+2$; 
make the two needed crossing changes to reach $\#_n(7_1\#\overline{7_1}$). 
But $u(L)=u(\overline{L})\geq 3n-1$, since anything smaller would 
allow us to unknot $\#_n 7_1$ more quickly. The $7$'s in this construction can
also be replaced
with any sequence of odd numbers $\geq 7$, while preserving the large additivity
gap, by Corollary \ref{cor:torus}.
What other examples can we build?

Even more than this, is it, perhaps, possible
that the unknotting number of a connected sum can be lower than that of its summand?

\begin{question}
Does $u(K\#K^\prime)\geq\max\{u(K),u(K^\prime)\}$ for all knots $K$ and $K^\prime$ ?
\end{question}

The fact that we could bootstrap an initial, single, example of 
non-additivity to find infinitely may 
knots $K$ for which $u(K\#7_1)<u(K)+u(7_1)$ has shown that 
failure of additivity is far from an isolated phenomenon.
This raises the question of whether or not every knot $K$ can be paired with a `poison' 
knot $K^\prime$ so that $u(K\#K^\prime)<u(K)+u(K^\prime)$. So we
ask:

\begin{question}
Does there exist a non-trivial knot $K$ for which 
$u(K\#K^\prime)=u(K)+u(K^\prime)$ for every knot $K^\prime$ ?
\end{question}

While we now know that 
$u(7_1\#\overline{7_1})<u(7_1)+u(\overline{7_1})$, and we can build
infinitely many examples of this phenomenon from this, 
all of our examples currently rely on
a single sequence of $5$ crossing changes that unknot $7_1\#\overline{7_1}$~.
We would certainly
like to understand why this sequence occurred! What, if
anything, is special about the knot $7_1$
that might make this happen? Understanding this unknotting sequence might help us 
understand how to find other methods to unknot knots more 
`efficiently'.

The intermediate
knot $K_1=K14a18636$, which played a central role in establishing the 
non-additivity of unknotting
number, has, from our work in \cite{bh21}, a weak BJ-unknotting 
number of 4, yet the unknotting sequence above shows that it has unknotting number 
at most 3. So the knot $K_1$ is, in fact, a counterexample to the weak
Bernhard-Jablan conjecture \cite{bh21}. That is, no crossing change in 
its (unique) minimal crossing diagram lowers the weak BJ-unknotting number.
This may point the way to building more examples
of unexpectedly low unknotting number, by trying to find a 
BJ-unknotting counterexample in an unknotting sequence. We note that 
there are (at least) six other 
knots that the authors have identified with weak BJ-unknotting 4 and 
unknotting number $\leq 3$: they are the knots $K13n1669$, $K14a2644$, 
$K14a11300$, $K14a11565$, $K14a12383$, and $K14a13809$. Each
is crossing-adjacent to a 15-crossing
diagram for a knot $K$ with $u(K)\leq 2$
(the corresponding unknotting number (at most) 2 knots are, respectively, 
$K15n9318$, $K15n9318$, $K15n35674$, $K15n82797$, $K15n4866$, and $K15n129095$),
and so each of these knots has unknotting number 
at most $3$. These could, like $K14a18636$, be useful targets
for building other unexpectedly short unknotting sequences for connected sums.

Finally, we end with the question:

\begin{question}
What is the actual unknotting number of $7_1\#\overline{7_1}$ ?
\end{question}

\noindent We now know that the unknotting number of $7_1\#\overline{7_1}$ is at most
 $5$, and, by Scharlemann's result \cite{scha85}, it is at least 2. Could it actually
be lower than $5$?

%%%%%%%%%%%%%%%%%%%%%%%%%%%%%%%%%%%%%%%%%%%%%%%%%%%%%%%%%%%%%%%%%%%%%%%%%%%%
%%%%%%%%%%%%%%%%%%%%%%%%%%%%%%%%%%%%%%%%%%%%%%%%%%%%%%%%%%%%%%%%%%%%%%%%%%%%

\section{Verification code}\label{sec:verify}

%%%%%%%%%%%%%%%%%%%%%%%%%%%%%%%%%%%%%%%%%%%%%%%%%%%%%%%%%%%%%%%%%%%%%%%%%%%%
%%%%%%%%%%%%%%%%%%%%%%%%%%%%%%%%%%%%%%%%%%%%%%%%%%%%%%%%%%%%%%%%%%%%%%%%%%%%

For the convenience of the reader, we give the code (written for Sage; 
remove the occurances of `snappy' to run in SnapPy) which will verify 
all of the identifications made for our example. We also include the output
from this code, as carried out in Sage. The reader's output may vary slightly,
due to some randomness built into some of SnapPy's routines. Note: When cutting 
and pasting from the code below, the quotation marks may need to be retyped
by the reader before pasting into SnapPy/Sage.

\medskip

\begin{verbatim}

import snappy

brd=[1,-4,2,3,3,3,2,3,2,2,4,-3,-3,-3,-3,-1,-3,-2,-3,-3]
print('First braid word (for L) is')
print(str(brd))
L=snappy.Link(braid_closure=brd)  ## building the connected sum
L.simplify('global')
print('Summands of L, and their number of crossings:')
PP=K.deconnect_sum()      ## this yields 7_1 and mirror(7_1)
print(str(PP))
print('Attempting to identify summands:')
print(str(PP[0].exterior().identify()))   ## these may be "[]"; SnapPy does 
print(str(PP[1].exterior().identify()))   ## not always succeed in 'recognizing' 
                              ## non-hyperbolic knots.
print('First summand knot group is ')
print(str(PP[0].exterior().fundamental_group()))
print('Second summand knot group is ')
print(str(PP[1].exterior().fundamental_group()))
## both of these will yield <a,b|aaaaaaabb> (or <a,b|aabbbbbbb>); 
## 7_1 is the only knot with this group.

brd2=brd[:]    ## implementing the 2 crossing changes
brd2[0]=-brd2[0]
brd2[1]=-brd2[1]
LA=snappy.Link(braid_closure=brd2)  ## this is K14a18636
print('Second braid word (for LA) is ')
print(str(brd2))
LA.simplify('global')
print('LA is '+str(LA.exterior().identify()))

DTC=[4,-16,24,26,18,20,28,22,-2,10,12,30,6,8,14]
DTCB=DTC[:]
print('DT code for KA is ')
print(str(DTC))
DTCB[0]=-DTCB[0]   ## change one crossing
print('DT code for KB is ')
print(str(DTCB))
KA=snappy.Link('DT:'+str(DTC))     ## this is K14a18636 
KB=snappy.Link('DT:'+str(DTCB))    ## this is K15n81556
print('KA is '+str(KA.exterior().identify()))
print('KB is '+str(KB.exterior().identify()))

print('SnapPy check: LA and KA are the same knot? '
+str(LA.exterior().is_isometric_to(KA.exterior())))

DTD=[4,12,-24,14,18,2,20,26,8,10,-28,-30,16,-6,-22]
print('DT code for KC is ')
print(str(DTD))
DTDB=DTD[:]
DTDB[6]=-DTDB[6]   ## change one crossing
print('DT code for KD is ')
print(str(DTDB))
KC=snappy.Link('DT:'+str(DTD))     ## this is K15n81556
KD=snappy.Link('DT:'+str(DTDB))    ## this is K12n412
print('KC is '+str(KC.exterior().identify()))
print('KD is '+str(KD.exterior().identify()))

print('SnapPy check: KB and KC are the same knot? '
+str(KB.exterior().is_isometric_to(KC.exterior())))

DTE=DTDB[:]
DTE[13]=-DTE[13]
print('DT code for KE is ')
print(str(DTE))
KE=snappy.Link('DT:'+str(DTE))    ## this is the unknot
print('KE is '+str(KE.exterior().identify()))

print('Knot group for KE:')
print(str(KE.exterior().fundamental_group()))
## only the unknot has cyclic knot group

\end{verbatim}

\noindent Sample output from this code:

\begin{verbatim}

>  First braid word (for L) is
>  [1, -4, 2, 3, 3, 3, 2, 3, 2, 2, 4, -3, -3, -3, -3, -1, -3, -2, -3, -3]
>  True
>  Summands of L, and their number of crossings:
>  [<Link: 1 comp; 7 cross>, <Link: 1 comp; 7 cross>]
>  Attempting to identify summands:
>  [7_1(0,0), K7a7(0,0)]
>  [7_1(0,0), K7a7(0,0)]
>  First summand knot group is
>  Generators:
>     a,b
>  Relators:
>     aabbbbbbb
>  Second summand knot group is
>  Generators:
>     a,b
>  Relators:
>     aabbbbbbb
>  Second braid word (for LA) is
>  [-1, 4, 2, 3, 3, 3, 2, 3, 2, 2, 4, -3, -3, -3, -3, -1, -3, -2, -3, -3]
>  True
>  LA is [K14a18636(0,0)]
>  DT code for KA is
>  [4, -16, 24, 26, 18, 20, 28, 22, -2, 10, 12, 30, 6, 8, 14]
>  DT code for KB is
>  [-4, -16, 24, 26, 18, 20, 28, 22, -2, 10, 12, 30, 6, 8, 14]
>  KA is [K14a18636(0,0)]
>  KB is [K15n81556(0,0)]
>  SnapPy check: LA and KA are the same knot? True
>  DT code for KC is
>  [4, 12, -24, 14, 18, 2, 20, 26, 8, 10, -28, -30, 16, -6, -22]
>  DT code for KD is
>  [4, 12, -24, 14, 18, 2, -20, 26, 8, 10, -28, -30, 16, -6, -22]
>  KC is [K15n81556(0,0)]
>  KD is [K12n412(0,0)]
>  SnapPy check: KB and KC are the same knot? True
>  DT code for KE is
>  [4, 12, -24, 14, 18, 2, -20, 26, 8, 10, -28, -30, 16, 6, -22]
>  KE is []
>  Knot group for KE:
>  Generators:
>     a
>  Relators:
>  

\end{verbatim}

%%%%%%%%%%%%%%%%%%%%%%%%%%%%%%%%%%%%%%%%%%%%%%%%%%%%%%%%%%%%%%%%%%%%%%%%%%%%
%%%%%%%%%%%%%%%%%%%%%%%%%%%%%%%%%%%%%%%%%%%%%%%%%%%%%%%%%%%%%%%%%%%%%%%%%%%%

\section{Acknowledgements}\label{sec:ackn}

%%%%%%%%%%%%%%%%%%%%%%%%%%%%%%%%%%%%%%%%%%%%%%%%%%%%%%%%%%%%%%%%%%%%%%%%%%%%
%%%%%%%%%%%%%%%%%%%%%%%%%%%%%%%%%%%%%%%%%%%%%%%%%%%%%%%%%%%%%%%%%%%%%%%%%%%%

The first author acknowledges support by a grant from
the Simons Foundation (Collaboration Grant number 525802).
The second author acknowledges support by a grant from
the Simons Foundation (Collaboration Grant number 581433).
The authors also acknowledge the support of
the Holland Computing Center at the University of 
Nebraska, which provided computing facilities on 
which the some of the work toward this project was carried out.

%%%%%%%%%%%%%%%%%%%%%%%%%%%%%%%%%%%%%%%%%%%%%%%%%%%%%%%%%%%%%%%%%%%%%%%%%%%%
%%%%%%%%%%%%%%%%%%%%%%%%%%%%%%%%%%%%%%%%%%%%%%%%%%%%%%%%%%%%%%%%%%%%%%%%%%%%
%%%%%%%%%%%%%%%%%%%%%%%%%%%%%%%%%%%%%%%%%%%%%%%%%%%%%%%%%%%%%%%%%%%%%%%%%%%%

%%%%%  THE BIBLIOGRAPHY

%%%%%%%%%%%%%%%%%%%%%%%%%%%%%%%%%%%%%%%%%%%%%%%%%%%%%%%%%%%%%%%%%%%%%%%%%%%%
%%%%%%%%%%%%%%%%%%%%%%%%%%%%%%%%%%%%%%%%%%%%%%%%%%%%%%%%%%%%%%%%%%%%%%%%%%%%


\begin{thebibliography}{99}


\bibitem{abe25}
T. Abe, personal communication (2025).

\bibitem{ahh12}
T. Abe, R. Hanaki and R. Higa
\emph{The uknotting number and band-unknotting number of a knot},
Osaka J. Math.
{\bf 49}
(2012)
523--550.

\bibitem{alex28}
J.W. Alexander,
\emph{Topological invariants of knots and links}, 
Trans. Amer. Math. Soc.
{\bf 30}
(1928)
275--306.

\bibitem{abdejltz24}
T. Applebaum, S. Blackwell, A. Davies, T. Edlich, A. Juh\'{a}sz, M. Lackenby, N. Toma\v{s}ev, and D. Zheng,
\emph{The unknotting number, hard unknot diagrams, and reinforcement learning}, 
Exp. Math. {\bf 35} (2025) 19pp.

\bibitem{ba05}
S. Baader,
\emph{Slice and Gordian numbers of track knots},
Osaka J. Math. 
{\bf 42} (2005)
257-–271.

\bibitem{ba06}
S. Baader,
\emph{Note on crossing changes}, 
Quarterly. J. Math. 
{\bf 57} 
(2006)
139--142.

\bibitem{ba10}
S. Baader,
\emph{Unknotting sequences for torus knots},
Math. Proc. Cam. Phil. Soc. 
{\bf 148}
(2010)
111--116.

\bibitem{ba25}
S. Baader, 
personal communication
(2025).

\bibitem{bm90}
J. S. Birman and W. Menasco, 
\emph{Studying links via closed braids IV: Composite links and split links}, 
Invent. Math. 
{\bf 102} 
(1990)
115–-139.

\bibitem{bl16}
M. Borodzik and C Livingston,
\emph{Semigroups, d -invariants and deformations of cuspidal singular points of plane curves},
J. London Math. Soc. 
{\bf 93} 
(2016)
439--463.

\bibitem{bh21}
M. Brittenham and S. Hermiller,
\emph{A counterexample to the Bernhard–Jablan Unknotting Conjecture},
Exp. Math. {\bf 30} (2021), no. 4, 547--556.

\bibitem{bh22}
M. Brittenham and S. Hermiller,
\emph{The smooth 4-genus of (the rest of) the prime knots through 12 crossings},
J. Knot Theory Ram. 
{\bf 31}
(2022)
Paper No. 2250081, 16 pp.

\bibitem{coli}
T. Cochran and W.B.R. Lickorish,
{\em Unknotting information from 4-manifolds}, 
Trans. Amer. Math. Soc. 
{\bf 297} 
(1986) 
125–-142. 

\bibitem{snappy}
M. Culler, N. Dunfield, M. Goerner, and J. Weeks, 
{\em SnapPy, a computer program for studying the geometry and topology of 3-manifolds}, 
http://snappy.computop.org .

\bibitem{dt83}
C. Dowker and M. Thistlethwaite,
\emph{Classification of knot projections},
Topology Appl. 
{\bf 16} 
(1983)
19--31.

\bibitem{fel14}
P. Feller, 
\emph{Gordian adjacency for torus knots}, 
Algebr. Geom. Topol. 
{\bf 14} 
(2014) 
769--793.

\bibitem{homfly}
P. Freyd, D. Yetter, J. Hoste, W.B.R. Lickorish, K. Millett, and A. Ocneanu,
\emph{A new polynomial invariant of knots and links},
Bull. Amer. Math. Soc.
{\bf 12} 
(1985)
239--246.

\bibitem{gl89}
C. McA. Gordon and J. Luecke, 
\emph{Knots are determined by their complements},
J. Amer. Math Soc.
{\bf 2}
(1989)
371--418.

\bibitem{plans77}
J.C. Hausmann (ed.),
\emph{Knot theory},
Proceedings of a Seminar held in Plans-sur-Bex, 1977,
Lecture Notes in Math.
{\bf 685}
(1978)
311 pp.

\bibitem{jones85}
V. Jones,
\emph{A polynomial invariant for knots via von Neumann algebras},
Bull. Amer. Math. Soc.
{\bf 12} 
(1985)
103--111.

\bibitem{kan10}
T. Kanenobu,
\emph{Upper bound for the alternation number of a torus knot},
Top. Appl.
{\bf 157}
(2010)
302--318.

\bibitem{kawa}
T. Kawamura,
{\em The unknotting numbers of $10_{139}$ and $10_{152}$ are $4$}, 
Osaka J. Math. 
{\bf 35} 
(1998)
539–-546.

\bibitem{ka96}
A. Kawauchi,
\emph{A survey of knot theory},
{Birkh\"auser Verlag, Basel},
{1996}.

\bibitem{kir93}
R. Kirby, 
\emph{Problems in low-dimensional topology}, 
in \emph{Geometric topology}, AMS/IP Studies in Advanced Mathematics, 2.2, 
American Mathematical Society, Providence, RI, 
(1997)

\bibitem{krmr}
P. Kronheimer and T. Mrowka,
{\em Gauge theory for embedded surfaces. I}, 
Topology 
{\bf 32} 
(1993)
773-–826. 

\bibitem{li02}
C. Livingston, 
\emph{The slicing number of a knot}, 
Algebr. Geom. Topol. {\bf 2} (2002) 1051–-1060.

\bibitem{knotinfo}
C. Livingston and A. Moore, \emph{KnotInfo: Table of Knot Invariants}, \par
{\tt knotinfo.org} (6/15/2025).

\bibitem{mc17}
D. McCoy,
{\em Alternating knots with unknotting number one},
Adv. Math. 
{\bf 305} 
(2017) 
757–-802.

\bibitem{mu65}
K. Murasugi, 
\emph{On a certain numerical invariant of link types}, 
Trans. Amer. Math. Soc. 
{\bf 117} 
(1965) 
387–-422. 

\bibitem{nn08}
S. Nelson and W. Neumann,
\emph{The 2-generalized knot group determines the knot},
Commun. Contemp. Math. 
{\bf 10} (2008)
843–-847.

\bibitem{ow08}
B. Owens, 
\emph{Unknotting information from Heegaard Floer homology},
Adv. in Mathematics 
{\bf 217}
(2008)
2353-–2376.

\bibitem{ow10}
B. Owens,
{\em On slicing invariants of knots},
Trans. Amer. Math. Soc. 
{\bf 362} 
(2010)
3095-–3106.

\bibitem{oza10}
M. Ozawa, 
\emph{Ascending number of knots and links},
J. Knot Theory Ram.
{\bf 19}
(2010)
15--25.

\bibitem{os03}
P. Ozsv\'{a}th and Z. Szab\'{o},
\emph{Knot Floer homology and the four-ball genus},
Geom. Topol. 
{\bf 7}
(2003)
615--639.

\bibitem{ozsz}
P. Ozsv\'{a}th and Z. Szab\'{o},
{\em Knots with unknotting number one and Heegaard Floer homology},
Topology 
{\bf 44} 
(2005)
705-–745.

\bibitem{ru93}
L. Rudolph,
\emph{Quasipositivity as an obstruction to sliceness},
Bull. Amer. Math. Soc. 
{\bf 29} 
(1993)
51–-59. 

\bibitem{sage}
The Sage Developers,
\emph{SageMath, the Sage Mathematics Software System (Version 9.2)},
{2025}, 
{\tt https://www.sagemath.org}.

\bibitem{scha85}
M. Scharlemann,
\emph{Unknotting number one knots are prime},
Invent. Math. 
{\bf 82} 
(1985)
37-–55.

\bibitem{schu53}
 H. Schubert, 
\emph{Knoten und Vollringe}, 
Acta Math. 
{\bf 90}
(1953)
131--286.

\bibitem{schu54}
H. Schubert, 
\emph{Uber eine numerische Knoteninvariante}, 
Math. Zeit.
{\bf 61} 
(1954)
245--288.

\bibitem{schu61}
H. Schubert, \emph{Bestimmung der Primfaktorzerlegung von Verkettungen}, 
Math. Zeit. 
{\bf 76} 
(1961)
116--148.

\bibitem{knotplot}
R. Scharein, 
\emph{Knotplot, a program for visualizing and interacting with 3-D and 4-D knots},
{\tt https://knotplot.com/}.

\bibitem{sima15}
V. Siwach and P. Madeti,
\emph{An unknotting sequence for torus knots},
Topology Appl. 
{\bf 196} 
(2015)
668--674.

\bibitem{sto02} 
A. Stoimenow,
\emph{The granny and the square tangle and the unknotting number},
Topology Appl. 
{\bf 117} 
(2002) 
59--75.

\bibitem{stoi}
A. Stoimenow,
\emph{Knot data tables}, available at 

{\tt https://stoimenov.net/stoimeno/homepage/ptab/} .

\bibitem{tan09}
K. Taniyama,
\emph{Unknotting numbers of diagrams of a given nontrivial knot are unbounded},
J. Knot Theory Ram.
{\bf 18}
(2009)
1049--1063.

\bibitem{vog90}
P. Vogel, 
\emph{Representation of links by braids: a new algorithm}, 
Comment. Math. Helv. 
{\bf 65} 
(1990)
104--113.

\bibitem{wa68}
F. Waldhausen,
\emph{On irreducible 3-manifolds which are sufficiently large},
Annals of Math.
{\bf 87} (1968) 
56–-88.

\bibitem{wz25}
C. Wang and Y. Zhang, 
\emph{A remark on the counterexample to the unknotting number conjecture},
preprint, 2 pp. {\tt arXiv:2507.14265}

\bibitem{wendt37}
H. Wendt,
\emph{Die gordische Aufl\"osung von Knoten},
Math. Zeit. 
{\bf 42} 
(1937)
680--696.

\bibitem{wh87}
W. Whitten,
\emph{Knot complements and groups},
Topology
{\bf 26}
(1987)
41--44.

\bibitem{ya08}
Z. Yang,
\emph{Unknotting number of the connected sum of $n$ identical knots},
J. Knot Theory Ram.
{\bf 17} 
(2008)
253--255.

\end{thebibliography}
\end{document}